\documentstyle[12pt]{article}
 \newtheorem{definition}{\sc Definition}[section]
 \newcounter{ex}
 \newenvironment{example}
   {{\sc Example \stepcounter{ex}\arabic{ex}\ :\ }}{}

\author{S.S.Moskaliuk, A.T. Vlassov  \\
   {\it   Bogolyubov Institute for Theoretical Physics} \\
      {\it 14b, Metrolohichna str., Kiev, UA-252143, Ukraine}\\
      {\it   E-mail: mss@gluk.apc.org }
}
\title{Double Categories in  Mathematical Physics}
\date{June 2, 2008}
\begin{document}
\maketitle


\begin{abstract}
Expansion of the categorical point of view on many areas of the mathematics
and mathematical physics will cause to deeper understanding of genuine
features of these problems. New applications of categorical methods are
connected with new additional structures on categories. One of such
structures, the double category, is considered in this article.
The double category structure is defined as generalization of the bicategory
structure. It is shown that double categories exist in the topological and
ordinary quantum field theories, and for dynamical systems with inputs and
outputs. Morphisms of all these double categories are not maps of sets.
\end{abstract}

\section{Double Category}

\begin{definition}
A double category $D$ consists of the following:

(1) A category $D_0$ of objects $Obj(D_0)$ and morphisms $Mor(D_0)$ of $0$%
-level.

(2) A category $D_1$ of objects $Obj(D_1)$ of $1$-level and morphisms $%
Mor(D_1)$ of $2$-level.

(3) Two functors $d,r:D_1\overrightarrow{\to }D_0.$

(4) A composition functor
$$
\ast :D_1\times _{D_0}D_1\to D_1
$$
where the bundle product is defined by commutative diagram
$$
\begin{array}{ccc}
D_1\times _{D_0}D_1 & \mathop{\rightarrow}\limits^{\pi _2} & D_1 \\
\pi _1\downarrow \quad &  & \quad \downarrow d \\
D_1 & \mathop{\rightarrow}\limits^{r} & D_0
\end{array}
$$

(5) A unit functor $ID:D_0\to D_1$, which is a section of $d,r$.

The above data is subject to {\bf Associativity Axiom} and {\bf Unit Axiom.}
If both of them are fulfilled only up to equivalence then the double
category is called a {\bf weak} double category, if they are fulfilled
strictly then it is a {\bf strong} double category.
\end{definition}

Here we see that for two objects $A,B\in Obj(D_0)$ there are $0$-level
morphisms $D_0(A,B)$ which we note by ordinary arrows $f:A\to B,$ and $1$%
-level morphisms $D_{(1)}(A,B)$ which we note by the arrows $\xi
:A\Rightarrow B,$ for $A=d(\xi )$ and $B=r(\xi )$. So with a $2$-level
morphism $\alpha :\xi \to \xi ^{\prime }$, where $\xi :A\Rightarrow B$ and $%
\xi ^{\prime }:A^{\prime }\Rightarrow B^{\prime }$ we can associate the
following diagram
$$
\begin{array}{ccccc}
A & \buildrel{\xi }\over{\Rightarrow } & B &  & \xi \\
d(\alpha )\downarrow \qquad &  & \qquad \downarrow r(\alpha ) & \qquad
\longmapsto \qquad & \quad \downarrow \alpha \\
A^{\prime } & \buildrel{\xi ^{\prime }}\over{\Rightarrow } & B^{\prime } &  & \xi
^{\prime }
\end{array}
$$

and arrow $\alpha :d(\alpha )\Rightarrow r(\alpha )$

The composition on 2-level associated with the diagram
$$
\begin{array}{ccccc}
A & \buildrel{\xi }\over{\Rightarrow } & B &  & \xi \\
d(\alpha )\downarrow \qquad &  & \qquad \downarrow r(\alpha ) &  & \quad
\downarrow \alpha \\
A^{\prime } & \buildrel{\xi ^{\prime }}\over{\Rightarrow } & B^{\prime } & \qquad
\longmapsto \qquad & \xi ^{\prime } \\
d(\alpha ^{\prime })\downarrow \qquad \; &  & \qquad \;\downarrow r(\alpha
^{\prime }) &  & \quad \;\downarrow \alpha ^{\prime } \\
A^{\prime \prime } & \buildrel{\xi ^{\prime \prime }}\over{\Rightarrow } &
B^{\prime \prime } &  & \xi ^{\prime \prime }
\end{array}
$$

Now we can define for double categories {\bf double (category) functors} and
their {\bf morphisms},{\bf double subcategories} , the category $DCat$ of
double categories, {\bf equivalence} of double categories, {\bf dual double
categories} (changed direction off 1-level morphisms, i.e. $d,r$ are
transposed), and so on.

\section{Examples of Double Categories}

Examples considered bellow show that double categories are sufficiently
natural for mathematics.

\begin{example}
Bicategories (\cite{Gab}) is the partial case of double category $D$
when the category $D_0$ is trivial, i.e. has only identical morphisms and
composition of 1-level and 2-level morphisms are associative.
\end{example}

\begin{example}
For each category $C$ we have the canonical double category $Morph(C)$ of
morphisms. Let $C$ be a category, $T$ be the diagram $\bullet \to \bullet ,$
$TC$ be the category of diagrams in $C$ of type $T$, let $D_0=C$ and $%
D_1=TC. $ The functor $d$ maps the diagram $f:A\to B$ into the object $A,$
the functor $r$ maps this diagram into the object $B,$ and so on. It is easy
to see that we get a double category $D$ which is noted by $Morph(C)$. Here $%
Obj(D_1)=Mor(D_0)$ , a 2-level morphism $f\Rightarrow g$ is a pair $(u,v)$
of morphisms $u,v\in Mor(C)$ from the commutative diagram
$$
\begin{array}{ccc}
A & \buildrel{u}\over{\rightarrow } & A^{\prime } \\
f\downarrow ~ &  & ~\downarrow f^{\prime } \\
B & \buildrel{v}\over{\rightarrow } & B^{\prime }
\end{array}
$$
with usual composition.
\end{example}

\begin{example}
Let $C$ be a category with bundle products, i.e. for all morphisms $u,v$ to $%
Y$ the universal square
$$
\begin{array}{ccc}
X\times _ZY & \longrightarrow & Y \\
\downarrow &  & \quad \downarrow v \\
X & \buildrel{u}\over{\longrightarrow } & Z
\end{array}
$$
exists. And let $T$ be the following diagram
$$
\bullet \leftarrow \bullet \to \bullet ,
$$
$TC$ be the category of diagrams in $C$ of type $T$. Now we define the
double category $D$ with $D_0=D$ and $D_1=TC.$ Two functors
$$
d,r:TC\to C,
$$
where the functor $d$ maps the diagram $A\leftarrow M\rightarrow B$ into the
object $A,$ the functor $r$ maps this diagram into the object $B.$ The
composition: for two $1$-level morphisms $\xi =(A\buildrel{\pi }\over{\leftarrow }M%
\buildrel{f}\over{\rightarrow }B):A\Rightarrow B$ and $\xi ^{\prime }=
(B\buildrel{\pi ^{\prime }}\over
{\leftarrow }M^{\prime ^{\prime }}\buildrel{f^{\prime }}\over{
\rightarrow }C):B\Rightarrow C$ we define their composition $\xi ^{\prime
}\circ \xi =(A\buildrel{\pi \circ \pi _1}\over{\leftarrow }M\times _BM^{\prime }
\buildrel{f\circ \pi _2}\over{\rightarrow }C)$ where the bundle product is defined
by the universal diagram
$$
\begin{array}{ccc}
{{M\times _BM^{\prime }}} & \mathop{\rightarrow}\limits^{\pi _2} & {M^{\prime }} \\
\pi _1\downarrow \quad &  & \quad \downarrow \pi ^{\prime } \\
M & \buildrel{f}\over{\rightarrow } & B
\end{array}
$$
A 2-level morphism is a triple $\alpha =(u,v,w):\xi \rightarrow \xi ^{\prime
}$ from the following commutative diagram
$$
\begin{array}{ccccc}
{M} & \buildrel{f}\over{\rightarrow } & B &  &  \\
{\buildrel{}\over{\pi \downarrow \ }} & \ \searrow v &  & \ \searrow w
&  \\
{\buildrel{}\over{A}} &  & M^{\prime } & \buildrel{f^{\prime }}\over{
\rightarrow } & B^{\prime } \\
& {\buildrel{}\over{\ \searrow u}} & \pi ^{\prime }\downarrow \ \  &  &
\\
&  & {{A^{\prime }}} &  &
\end{array}
$$
with the evident composition.
\end{example}

\begin{example}
 Let $k$ be a ring, $Alg_k$ be the category of unital $k$
-algebras. We define the double category with $D_0=Alg_k,$ $Obj(D_1)$ as a
set of bimodules. So that for a left $A$- and right $B$-module (= $A\otimes
_kB^{\circ }$-module) $M$ we suppose $d(M)=A,$ $r(M)=B.$ Composition of two
morphisms $A\buildrel{M}\over{\Rightarrow }B
\buildrel{N}\over{\Rightarrow }C\quad $
then defined as $\quad N*M=M\otimes _BN.$

For $A\buildrel{M}\over{\Rightarrow }B$ and
$C\buildrel{N}\over{\Rightarrow }D$ the
set $D_1(M,N)$ is the set of triples $(u,v,w),$ where $u:A\rightarrow
C,v:B\rightarrow D,w:M\rightarrow N$ such that $w(amb)=u(a)\cdot w(m)\cdot
w(b).$ Unit for $A\in Obj(A$l$g_k)$ is $A\otimes _kA^{\circ }$-module $A$.
Subcategories, see \cite{Lod}. There is the natural inclusion of double
categories
$$
Morph(A\mbox{l}g_k)\hookrightarrow ALG_k.
$$

Here 1-level morphism $N:A\Rightarrow B$ we can consider as an algebra
homomorphism $A\otimes _kB^{\circ }\rightarrow End_k(N)$ (what is a selected
element $n_0$?). Analogical formulas appear in the next more general
situation.
\end{example}

\begin{example}
For a multiplicative (tensor) category $(C,\otimes ,U,u)$ (see \cite{Del}).
Then we have the double category with $D_1=C,$ and $D_0=(*,*)$, trivial
category with one object and one morphism. The composition is
$$
D_1\times _{D_0}D_1=C\times C\buildrel{\otimes }\over{\rightarrow }C.
$$
Other double category is $D$ with $D_0=C$ and $D_1$ such that
$$
Obj(D_1)=\{(X,x)|A,B,X\in Obj(C),\quad \/x:X\otimes A\rightarrow B\}.
$$
So, we write $\xi =(X,x):A\Rightarrow B$ and for $\xi \in Obj(D_1)$ we
denote $\xi =(X_\xi ,x_\xi ),$ $d(\xi )=A_\xi ,\quad r(\xi )=B_\xi .$
2-level morphisms
$$
D_1(\xi ,\xi ^{\prime })=\{(f_1,f_2,f_3)\ |\ commutative\;diagram\;
\begin{array}{ccc}
X\otimes A & \buildrel{x}\over{\longrightarrow } & B \\
f_3\otimes f_1\downarrow \qquad ~\ \  &  & \quad \downarrow f_2^{\prime } \\
X^{\prime }\otimes A^{\prime } & \buildrel{x^{\prime }}\over{\longrightarrow } &
B^{\prime }
\end{array}
\}
$$
and $d(f_1,f_2,f_3)=f_1,\quad r(f_1,f_2,f_3)=f_2.$

Composition $D_1\times _{D_0}D_1\rightarrow D_1$ is defined so

for $A\buildrel{\xi }\over{\Rightarrow }B\buildrel{\xi ^{\prime }}
\over{\Rightarrow }%
B^{\prime }$ $\quad \xi \circ \xi ^{\prime }=(A,B^{\prime },X^{\prime
}X,x^{\prime \prime }),$ where $x^{\prime \prime }$ is the following
composition

$(X^{\prime }\otimes X)\otimes A\buildrel{\varphi _{X^{\prime },X,A}^{-1}}
\over{
\longrightarrow }X^{\prime }\otimes (X\otimes A)\buildrel{id_{X^{\prime
}}\otimes x}\over{\longrightarrow }X^{\prime }\otimes B
\buildrel{x^{\prime }}\over{\longrightarrow }B^{\prime }$
\end{example}

\section{Action of a double category}

Double categories are categorical variants of  usual monoids (and groups),
and thus we have the corresponding variant for their actions. Below the
definition of action of a double category $d,r:D_1\rightarrow D_0$ on
categories over $D_0$ is given. Thus we get an analog of 
group-theoretical methods in
categorical frames.

\begin{definition}
(Left) action of a double category $d,r:D_1\rightarrow D_0$ on a category $%
p:M\rightarrow D_0$ over $D_0$ is a functor $\varphi $ such that

(1)The diagram is commutative
$$
\begin{array}{ccc}
D_1\times _{D_0}M & \buildrel{\varphi }\over{\rightarrow } & M \\
&  &  \\
\quad \quad r\circ \pi _1 & \searrow & \quad \downarrow \;p \\
&  &  \\
&  & D_0
\end{array}
$$
where the bundle product $D_1\times _{D_0}M$ is defined by the diagram
$$
\begin{array}{ccc}
D_1\times _{D_0}M & \buildrel{\pi _2}\over{\rightarrow } & M \\
\pi _1\downarrow \quad &  & \quad \downarrow \;p \\
D_1 & \buildrel{r}\over{\rightarrow } & D_0
\end{array}
$$

(2)The diagram is commutative with precise to an isomorphism
$$
\begin{array}{ccccc}
(D_1\times _{D_0}D_1)\times _{D_0}M & \buildrel{\cong }\over{\rightarrow } &
D_1\times _{D_0}(D_1\times _{D_0}M) & \buildrel{id_{D_1}\times _{D_0}\varphi
}\over{\longrightarrow } & D_1\times _{D_0}M \\
\otimes \times _{D_0}id_M\downarrow \qquad \quad \qquad &  &  &  & \ \
\downarrow \varphi \\
D_1\times _{D_0}M &  & \buildrel{\varphi }\over{\longrightarrow } &  & M
\end{array}
$$
and there exists a functor isomorphism $\varphi $ such that
$$
\forall \;\xi ,\xi ^{\prime }\in Obj(D_1),\;m\in Obj(M_1)\;\quad \varphi
_{\xi ,\xi ^{\prime },m}:(\xi *\xi ^{\prime })*m\rightarrow \xi *(\xi
^{\prime }*m)
$$

(3) For the unit functor we have a functor isomorphism $\chi :\varphi \circ
(ID\times id_M)\widetilde{\longrightarrow }id_M$ or for objects
$$
\forall \;A\in Obj(D_0),\;m\in Obj(M_1)\quad \chi _{A,m}:ID_A*m\widetilde{%
\longrightarrow }m
$$
\end{definition}

So we have the map of pair of objects $\xi \in Obj(D_1),$ $m\in Obj(M)$ $(A%
\buildrel{\xi }\over{\Rightarrow }p(m),m)\mapsto \varphi (\xi ,m)$ such that $%
p(\varphi (\xi ,m))=A,$ and of morphisms $\alpha \in D_1(\xi ,\xi ^{\prime
}),u\in M(m,m^{\prime })$
$$
\begin{array}{ccccccc}
\xi & \qquad & A & \buildrel{\xi }\over{\Rightarrow } & p(m) &  &
 {\varphi (\xi ,m)} \\
\alpha \downarrow \quad &  & f=d(\alpha )\downarrow ~ &  & \quad \downarrow
r(\alpha )=p(u) & \quad \longmapsto \quad & \qquad \downarrow \varphi
(\alpha ,u) \\
\xi ^{\prime } &  & A^{\prime } & \buildrel{\xi ^{\prime }}
\over{\Rightarrow } &
p(m^{\prime }) &  & {\varphi (\xi ^{\prime },m^{\prime })}
\end{array}
$$
and here $p(\varphi (\alpha ,u))=f.$

The definition of a {\bf right action} is evident.

Invariant subcategories, ''orbits'', invariant properties.

An orbit of a subcategory is a subcategory.

\begin{example}
Each double category $D$ acts on itself one from left and from right by the
composition $*$.
\end{example}

\begin{example}
Let $SC$ be a subcategory of morphisms of the category $C$ such that for all
($\pi :M\rightarrow B)\in Obj(SC)$ and all $(f:B^{\prime }\rightarrow B)\in
Mor(C)$ there exists the universal square
$$
\begin{array}{cccc}
f^{*}M= & B^{\prime }\times _BM & \buildrel{f^{\prime }}\over{\longrightarrow } & M
\\
& \pi _1\downarrow \quad \; &  & \quad \downarrow \pi \\
& B^{\prime } & \buildrel{f}\over{\longrightarrow } & B
\end{array}
$$
and $\Pi :SC\rightarrow C$ is the projection on base, i.e. $\Pi :(\pi
:M\rightarrow B)\longmapsto B.$ Then we have the natural action of $Morph(C)$
on $SC$
$$
Morph(C)\times _CSC\longrightarrow SC
$$
which maps the pair $((f:B^{\prime }\rightarrow B),(\pi :M\rightarrow B))$
to $(\pi _1:f^{*}M\rightarrow B^{\prime }).$
\end{example}

\begin{example}
Let $mod_k$ be the category of left modules over $k$-algebras and
$\Pi :mod_k\rightarrow A$l$g_k$ be the natural projection, which an object $(A,M)$
maps to $A.$ There is natural action of $ALG_k$ on $mod_k$
$$
(ALG_k)_1\times _{(ALG_k)_0}m\mbox{o}d_k\rightarrow mod_k
$$
such that for $\xi =N:A\Rightarrow B$ and $B$-module $M$
$$
\xi ^{*}(B,M)=(A,N\otimes _BM).
$$
\end{example}

\begin{example}
Characteristic classes. Let a double category $G$ act on $p:M\rightarrow G_0$,
and there is a contravariant functor $H:G^{\circ }\rightarrow Morph(M)$. A
characteristic class of $m\in Obj(M)$ is $c(m)\in H_0(p(m)),$ such that for
all $\xi :p(m^{\prime })\Rightarrow  p(m)$ we have
$$
c(\xi ^{*}m)=H_1(\xi )(c(m)).
$$
\end{example}

\begin{example}
Equivariant functors. Let $M$ be a category of manifolds (topological or
smooth), $L$ be a category of locally trivial bundles over objects of $M$.
Then $Morph(M)$ acts on $L$. Let $G$ be a topological
group, $M_G$ be the category of $G$-manifolds, $P$ is the category of
principal bundles with structure group $G$ over objects of $M$. The functor
$$
P\times M_G\rightarrow L
$$
which maps $(\eta ,F)$ to the fiber bundle $\eta [F]$ with fiber $F.$ This
functor is equivariant relatively of action $Morph(M)$ on $P$ and $L.$
\end{example}

\begin{example}
Let $ISO({\cal C})$ be a sub double category of $Morph({\cal C})$ such that $%
(ISO({\cal C}))_0={\cal C}$ and $(ISO({\cal C}))_1$ is the full subcategory
of $(Morph({\cal C}))_1$ with
$$
Obj((ISO({\cal C}))_1)=\{f\in Mor({\cal C})\quad |\quad f\mbox{ is an
isomorphism }\}.
$$
For any forgetful functor $F:{\cal C}^{\prime }\rightarrow {\cal C}$ (see
the next section) the double category $ISO({\cal C})$ acts on ${\cal C}%
^{\prime }$ from left
$$
ISO({\cal C})\times _{{\cal C}}{\cal C}^{\prime }\rightarrow {\cal C}%
^{\prime }:(u:B\rightarrow F(C),C)\mapsto u^{*}C
$$
and $v^{*}(u^{*}C)\cong (u\circ v)^{*}C$.
\end{example}

\section{Cobordism and Double Categories}

Let $M_d$ be the category of oriented compact $d$-dimensional smooth
manifolds (with boundary) and piecewise smooth maps (the sense of the
condition we do not define more exactly here; this may be such continuous
maps $f:M\rightarrow Y$ that are smooth on a dense open subset $U_f\subset
M $ ) , let $CM_d$ be its subcategory of closed (with empty boundary)
manifolds and smooth maps, $CM_d\subset M_d$.

There are the following functors:

(1) Disjoint union
$$
\cup :M_d\times M_d\to M_d:(X,Y)\mapsto X\cup Y.
$$
(2) Changing of the orientation of manifolds on opposite
$$ (-):M_d\to M_d:X\mapsto -X.
$$
(3) Boundary operator
$$ \partial :M_{d+1}\to CM_d:X\mapsto \partial X. $$
(4) Multiplication on the unit segment $I=[0,1]$
$$ I\times \ .:CM_d\to M_{d+1}:X\mapsto I\times X. $$

Now we define a double category $\bf{C(d)}$ with

\begin{itemize}
\item[(1)]  $\bf{C(d)}_0=CM_d.$

\item[(2)]  1-level morphisms $\bf{C(d)}_{(1)}(X,X^{\prime })$ is a set
of pairs $(Y,f)$ where $Z$ is oriented compact $(d+1)$-dimensional smooth
manifold with the boundary $\partial Y$ and $f$ is an diffeomorphism
$$
f:(-X)\cup X^{\prime }\to \partial Y,
$$
where $\cup $ notes the disjoint union of $-X$ and $X^{\prime }.$ Thus we
write $(Y,f):X\Rightarrow X^{\prime }.$

\item[(3)]  The composition of $(Y,f):X\Rightarrow X^{\prime }$ and $%
(Y^{\prime },f^{\prime }):X^{\prime }\Rightarrow X^{\prime \prime }$ is the
morphism
$$ (Y\cup _{X^{\prime }}Y^{\prime },(f|_X)\cup (f^{\prime }|_{X^{\prime
}})):X\Rightarrow X^{\prime \prime },
$$
where $(Y\cup _{X^{\prime }}Y^{\prime })$ denotes the union $(Y\cup
Y^{\prime })$ after identification of each point $f(y)\in f(Y)$ with the
point $f^{\prime }(y)\in f^{\prime }(Y)$ for all $y\in Y$ and smoothing this
topological manifold.

\item[(4)]  The 1-level identical morphism $ID_X$ is $(X\times
[0;1],id_{(-X)\cup X}),$ because $\partial (X\times [0;1])=(-X)\cup X.$

\item[(5)]  2-level morphisms of $\bf{C(d)}_1(\xi ,\xi ^{\prime })$ from
$\xi =(Y,f:X^{\prime }\cup (-X)\to \partial Y):X\Rightarrow X^{\prime }$ to
$\xi ^{\prime }=(Y^{\prime },f^{\prime }:X^{\prime \prime }\cup (-X^{\prime
})\to \partial Y^{\prime }):X^{\prime }\Rightarrow X^{\prime \prime }$ are
such triples of smooth maps $(f_1,f_2,f_3)$ that the following diagram is
commutative

$$
\begin{array}{cccc}
{(-X)\cup X}^{\prime } & \buildrel{f}\over{\longrightarrow } & {\partial Y} &
\subset {Y} \\
{\quad \quad \;\quad \quad \downarrow f_1\cup f_2} &  &  & {{}{}{}{\quad
\downarrow f_3}} \\
{{(-X^{\prime })\cup X}}^{\prime \prime } & \buildrel{f^{\prime }}\over{
\longrightarrow } & {\partial Y}^{\prime } & \subset {Y}^{\prime }
\end{array}
$$
\end{itemize}

It easy to see that functors $\cup $ and $(-)$ may be expanded to double
category functors
$$
\begin{array}{c}
\cup :\bf{C(d)}\rightarrow \bf{C(d)}, \\
(-):\bf{C(d)}\rightarrow \bf{C(d)}^{\circ }
\end{array}
$$
and $(-)$ is an equivalence of the double categories.

{\bf Remark.}
It is interesting the appearance the following two formulas for 1-level
morphisms in algebras  and cobordisms
$$
f:A\otimes _kB^{\circ }\rightarrow End_k(N)\qquad f:(-X)\cup Y\to \partial Z,
$$
where we have correspondence between the functors
$$
\begin{array}{ccc}
(\_)^{\circ } & \quad \longleftrightarrow \quad & -(\_), \\
\otimes _k & \quad \longleftrightarrow \quad & \cup , \\
End_k & \quad \longleftrightarrow \quad & \partial .
\end{array}
$$

\section{Topological Quantum Field Theory}

Topological quantum field theory is a functor $Z$ from the category $CM(d)$
of $d$-dimensional manifolds to the category of $H$ of (usually Hermitian)
finite dimensional vector spaces and some axioms are satisfied (\cite
{Ati}). Really, the {\bf functor} $Z$ is a functor between double
categories.

Thus, topological quantum field theory in dimension $d$ is a functor
$$
Z:\bf{C(d)}\rightarrow Morph(H),
$$
between double categories such that:

(1) the disjoint union in $\bf{C(d)}$ go to the tensor product
$$
\cup \quad \mapsto \otimes ,
$$
where $(\_)^{*}:H\rightarrow H^{\circ }$ is dualization of vector spaces.

(2) changing of orientation in $\bf{C(d)}_0$ go to dualization
$$
(-)\mapsto (.)^{*}
$$

Thus, as consequence of double categorical functorial properties, we get

\begin{itemize}
\item[(1)]  for each compact closed oriented smooth d-dimensional manifold $%
X\in Obj(\bf{C(d)}_0)$ the value of the functor $Z(X)$ is a finite
dimensional vector space over the field $\bf{C}$ of the complex numbers
(usually with Hermitian metric),

\item[(2)]  for each $(Y,f):X\Rightarrow X^{\prime }$ from $Obj(\bf{C(d)%
}_1)$ the value of the functor $Z(Y,f)$ is a homomorphism $Z(X)\rightarrow
Z(X^{\prime })$ of (Hermitian) vector spaces$,$
\end{itemize}

and the following well known axioms of topological quantum field theory are
satisfied:

\begin{itemize}
\item[A(1)]  (involutivity) $Z(-X)=Z(X)^{*},$ where $-X$ denotes the
manifold with opposite orientation, and $*$ denotes the dual vector space.

\item[A(2)]  (multiplicativity) $Z(X\cup X^{\prime })=Z(X)\otimes
Z(X^{\prime }),$ where $\cup $ denotes disconnected union of manifolds.

\item[A(3)]  (associativity) For the composition $(Y^{\prime \prime
},f^{\prime \prime })=(Y,f)*(Y^{\prime },f^{\prime })$ of cobordisms must be
$$
Z(Y^{\prime \prime },f^{\prime \prime })=Z(Y^{\prime },f^{\prime })\circ
Z(Y,f)\in Hom_{\bf{C}}(Z(X),Z(X^{\prime \prime })).
$$
(Usually the identifications
$$
Z(X^{\prime }-X)\cong Z(X)^{*}\otimes Z(X^{\prime })\cong Hom_{\bf{C}%
}(Z(X),Z(X^{\prime }))
$$
allow us to identify $Z(Y,f)$ with the element $Z(Y,f)\in Z(\partial Y)$.

\item[A(4)]  For the initial object $\emptyset \in Obj(\bf{C(d)}%
_0)\qquad Z(\emptyset )=\bf{C}.$

\item[A(5)]  (trivial homotopy condition) $Z(X\times [0,1])=id_{Z(X)}.$
\end{itemize}

\subsection{Field Theory}

Here is sketch of categorical construction with double categories for the
ordinary field theory, where we deal with fiber bundle and equations for
their sections. Let us denote corresponding double category by $F(d)$where $%
(d+1)$ is the dimension of the space-time, and it is similar to the double
category $\bf{C(d)}.$

Objects of category $F(d)_0$ are $A=(\pi :V\rightarrow X,s\epsilon \Gamma
(V/X))$ where $X$ is an oriented compact closed $d$-dimensional manifold, $%
\pi :V\rightarrow X$ is a fiber bundle over $X$ of some defined type, $s$ is
a section of $\pi $ with some special properties. The definition of morphism
is evident but there are variants. There are functors $\cup $ and $-$ as in $%
\bf{C(d)}_0$.

Objects of $F(d)_1$ are defined the following way. Let $\pi :E\rightarrow Y$
be a fiber bundle over an oriented compact $(d+1)$-dimensional manifold $X$
with a boundary, $\cal{D}$ is a system of equations for sections of $\pi
,s$ is a solution of $\cal{D}$. On the category of these data we define
the boundary functor $\partial ,$ which maps an object $(\pi :E\rightarrow Y,%
\cal{D},s)$ to some fiber bundle $\widetilde{\partial }E\rightarrow
\partial Y$ with its section $\widetilde{\partial }s$. Here the construction
of $\widetilde{\partial }E$ depends of type of the system $\cal{D}$. An
object of $F(d)_1$ is $\xi =(A,A^{\prime },(E/Y,\cal{D},s),f)$ where $f$
is isomorphism
$$
f:(-A)\cup A^{\prime }\rightarrow \partial (E/Y,\cal{D},s).
$$
There exit the composition $\xi *\xi ^{\prime }.$ For the action integral $%
S[\xi ]$ we have $S[\xi *\xi ^{\prime }]=S[\xi ]+S[\xi ^{\prime }].$ The
system $\cal{D}$ is an Euler-Lagrange equation for action $S,$ and
boundary bundles and sections define unique gluing.

Thus it maybe we get categorical field theory.

\end{document}